\documentclass[notitlepage,12pt]{article} 
\usepackage{amssymb,amsmath,latexsym,epsfig,euscript,rotating,epic} 

\newlength{\jmr}
\newlength{\koi}
\newlength{\hwl}
\newlength{\khov}
\newlength{\bernd}
\newtheorem{lemma}{Lemma}

\newtheorem{thm}{Theorem}
\newtheorem{prop}{Proposition}
\newtheorem{dfn}{Definition}
\newtheorem{cor}{Corollary}
\newtheorem{rem}{Remark}	
 
\newtheorem{ex}{Example}

\newcommand{\pspa}{{\mathbf{PSPACE}}}

\newcommand{\am}{{\mathbf{AM}}}

\newcommand{\np}{{\mathbf{NP}}}

\newcommand{\conp}{{\mathbf{coNP}}}
\newcommand{\corp}{{\mathbf{coRP}}}
\newcommand{\rp}{{\mathbf{RP}}}
\newcommand{\bk}{{\bar{\kappa}}}
\newcommand{\kk}{{\kappa}}
\newcommand{\bp}{{\mathfrak{p}}}
\newcommand{\bq}{{\mathfrak{q}}}
\newcommand{\bpp}{{\mathbf{BPP}}}
\newcommand{\crap}{\pp^{\np^\np}}
\newcommand{\pp}{\mathbf{P}}
\newcommand{\hn}{\mathbf{HN}}

\newcommand{\expt}{{\mathbf{EXPTIME}}}

\newcommand{\eps}{\varepsilon}

\newcommand{\cO}{\mathcal{O}}

\newcommand{\gal}{\mathrm{Gal}}

\newcommand{\grh}{{\rm GRH }}

\newcommand{\thth}{^{\underline{\mathrm{th}}}}

\newcommand{\st}{ {\underline{ \mathrm{st} } }  }
\newcommand{\nd}{{\underline{\mathrm{nd}}}}

\newcommand{\Q}{\mathbb{Q}}
\newcommand{\R}{\mathbb{R}}
\newcommand{\C}{\mathbb{C}}
\newcommand{\N}{\mathbb{N}}
\newcommand{\Z}{\mathbb{Z}}
 
\newcommand{\cp}{\mathfrak{p}} 
\newcommand{\ca}{\mathfrak{a}}

\newcommand{\Rn}{\R^n}

\newcommand{\Cn}{\C^n}

\newcommand{\qed}{$\blacksquare$}
\newcommand{\dia}{$\diamond$}

\newcommand{\bO}{\mathbf{O}}
\newcommand{\vol}{\mathrm{Vol}}

\begin{document}

\title{\mbox{}\\ 
\vspace{-1in}
\mbox{}\hspace{-.4cm}\scalebox{.75}[1]{Dedekind Zeta Functions and the 
Complexity of Hilbert's Nullstellensatz} } 

% \subjclass{Primary 
% 11G25; % Arithmetic Algebraic Geometry; varieties over finite and 
%       % local ground fields
% Secondary 
% 11G35, % ...global ground fields
% 14D10, % algebraic geometry; families, fibrations; arithmetic ground fields
% 14G20. % algebraic geometry; arithmetic problems/diophantine geometry; 
%        % local ground fields 
% } 

% \author{
% \alignauthor 
\author{J.\ Maurice Rojas%\titlenote{
% \thanks{
%  This research was partially supported by 
%\\ 
\thanks{ 
% \affaddr{
Department of Mathematics, 
% }
% \affaddr{
Texas A\&M University, 
% } 
% \affaddr{
TAMU 3368, 
% } 
% \affaddr{
College Station, Texas 77843-3368, 
% }
% \affaddr{
USA. 
% }
% \affaddr{FAX: (979) 845-6028 }\\  
% \email{
e-mail: {\tt rojas@math.tamu.edu} , 
% \email{
web page: {\tt http://www.math.tamu.edu/\~{}rojas} \hfill . \hfill \mbox{}
This research was partially supported by 
NSF grant DMS-0211458 a grant from the Texas A\&M College of Science.}}

% \dedicatory{In memory of Ruth Ingrid Michler. }  

\maketitle 

\vspace{-.25in}
\begin{abstract} 
Let $\hn$ denote the problem of determining whether 
a system of multivariate polynomials with integer 
coefficients has a complex root. It has long been known that 
$\hn\!\in\!\pp \Longrightarrow \pp\!=\!\np$ and, thanks to recent work 
of Koiran, it is now known that the truth of the 
{\bf Generalized Riemann Hypothesis (GRH)} yields 
the implication $\hn\!\not\in\!\pp \Longrightarrow 
\pp\!\neq\!\np$. We show that the assumption of \grh in the latter implication 
can be replaced by either of two more plausible hypotheses from analytic number 
theory. The first is an effective short interval Prime Ideal Theorem with 
explicit dependence on the underlying field, while the second can be 
interpreted as a quantitative statement on the higher moments of the zeroes 
of Dedekind zeta functions. In particular, both assumptions can still hold 
even if \grh is false. We thus obtain a new application of Dedekind 
zero estimates to computational algebraic geometry. 
Along the way, we also apply recent explicit algebraic and analytic estimates, 
some due to Silberman and Sombra, which may be of independent interest. 
% note!: the small degree result needs some discriminant 
%        considerations NOT yet made in the Mitsui-Sokolovskii-Heath-Brown 
%        results... (June 13, 2002) 
% note!: the degree can be approximated to within a 
%        constant factor via my stockmeyer trick from 
%        jcs.tex !!! (i used deg>1 detection to solve 
%        transitivity detection) maybe I should write this in a 
%        separate paper? !!! 
\end{abstract} 

% note: x^(20/9) grows quickly enough 
%       for the K=Q version of the stride theorem! 

\section{Introduction and Main Results} 
\label{sec:intro} 
Deciding whether a system of multivariate polynomial equations  
has a solution is arguably the most fundamental problem of 
algebraic complexity. Over the ring $\Z/2\Z$ this problem 
defines the famous complexity class $\np$ and is  
the subject of much interest  in proof 
systems theory (see, e.g., \cite{bikpp}). Over the complex numbers this 
problem is a cornerstone of computational algebraic geometry and is 
the defining problem of the complexity class $\np_\C$ of 
the recent BSS model over $\C$ \cite{bcss}. 
Let $\hn$ denote the problem of determining whether
a multivariate system of polynomials with {\bf integer} 
coefficients has a {\bf complex} root. It has been known 
at least since the 1980's that $\hn\!\in\!\pspa$ (see, e.g., \cite{pspace}), 
and that the special case of two univariate polynomials is 
already $\np$-hard \cite{plaisted}. 

\begin{rem} 
Recall the following containments of complexity classes:\\ 
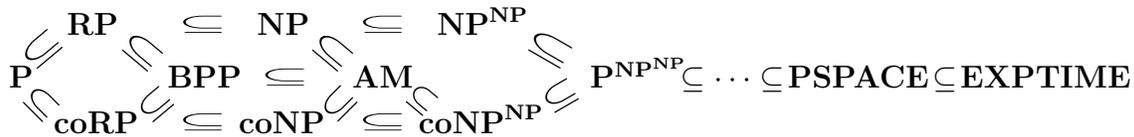
\begin{picture}(133,80)(0,10)
\put(0,35){$\pp$}
\put(8,42){\begin{rotate}{40}\scalebox{1.8}[1]{$\subseteq$}\end{rotate}}
\put(22,55){$\rp$} 
\put(5,32){\begin{rotate}{-45}\scalebox{1.5}[1]{$\subseteq$}\end{rotate}}
\put(17,18){$\corp$}
\put(42,52){\begin{rotate}{-40}\scalebox{1.8}[1]{$\subseteq$}\end{rotate}}
\put(52,22){\begin{rotate}{40}\scalebox{1.5}[1]{$\subseteq$}\end{rotate}}
\put(60,35){$\bpp$}
\put(64,57){\scalebox{2}[.7]{$\subseteq$}}
\put(94,55){$\np$} 
\put(114,52){\begin{rotate}{-40}\scalebox{1.8}[1]{$\subseteq$}\end{rotate}} 
\put(64,20){\scalebox{2}[.7]{$\subseteq$}}
\put(87,18){$\conp$} 
\put(132,20){\scalebox{2}[.7]{$\subseteq$}} 
\put(122,22){\begin{rotate}{45}\scalebox{1.5}[1]{$\subseteq$}\end{rotate}} 
\put(95,37){\scalebox{2}[.7]{$\subseteq$}} 
\put(130,35){$\am$} 
\put(150,32){\begin{rotate}{-40}\scalebox{1.2}[1]{$\subseteq$}\end{rotate}} 
\put(132,57){\scalebox{2}[.7]{$\subseteq$}} 
\put(162,55){$\np^{\np}$} 
\put(195,52){\begin{rotate}{-30}\scalebox{2}[1]{$\subseteq$}\end{rotate}} 
\put(155,18){$\conp^{\np}$} 
\put(204,26){\begin{rotate}{30}\scalebox{1.5}[1]{$\subseteq$}\end{rotate}} 
\put(220,35){$\crap$} 
\put(254,35){$\subseteq\cdots \subseteq\!\pspa\!\subseteq\!\expt$} 
\end{picture}\\
and that all these complexity classes consist of {\bf decision} problems, 
i.e., problems where the answer to any instance is either {\tt Yes} or 
{\tt No}. In particular, the properness of {\bf every} inclusion above 
is a major open problem, and even the inclusion $\pp\!\subseteq\!\pspa$ is 
not known to be proper \cite{zachos,lab,arith,papa}. 
% Informally, the 
% underlying computational model can be taken to be a lap-top computer 
% with flawless hardware and infinite memory, and the 
% notions of ``time'', ``work'', or ``number of steps'' can be taken as 
% the number of clock-cycles of 
% the underlying microprocessor needed to solve a given problem 
% instance.\footnote{The equivalence 
% of such a machine to the classical {\bf Turing machine} is simply a 
% restatement of the well-known RAM model of computation \cite{papa}. 
% The equivalence of a clock-cycle to a small number of ``bit operations'' 
% follows from the basic engineering principles of all modern 
% microprocessors. } 
The definitions of the aforementioned complexity classes are reviewed briefly in 
Section \ref{sec:back}. \dia 
\end{rem} 

Here we use analytic number theory to get a simple, surprisingly fast 
algorithm for $\hn$. Our algorithm is completely different from the 
usual Gr\"obner basis, resultant, or homotopy methods from computational 
algebraic geometry. Indeed, the latter methods result in a complexity 
bound of at best $\pspa$ for deciding the existence of a 
complex root (see, e.g., \cite{pspace,four}), while our algorithm yields a bound of 
$\am$ --- the so-called {\bf Arthur-Merlin} class \cite{lab} --- under either of two  
plausible number-theoretic assumptions: $\pmb{\text{GIPIT}}$ and $\pmb{\text{DZH}}$. 
In particular, 
an Arthur-Merlin bound for $\hn$ was known earlier only under the assumption 
of the {\bf Generalized Riemann Hypothesis (GRH)} \cite{hnam} (see the next section 
for its statement), while our new hypotheses can both continue to hold even 
if \grh is false. 

In essence, our approach begins by naturally associating a finite algebraic extension $K$ 
of $\Q$ to any polynomial system $F$. Then, if the prime ideals in the ring of integers 
$\cO_K$ of $K$ are sufficiently well behaved, we can use a randomized sampling trick to 
relate the existence of complex roots of $F$ to the existence of roots of 
a mod $p$ reduction of $F$. Our hypotheses GIPIT and DZH (short 
for {\bf Good Short Interval Prime Ideal Theorem} and {\bf Dedekind 
Zero Hypothesis}) are then simply the two most convenient statements encoding what sort 
of quantitative behavior we require from the prime ideals of $\cO_K$. To 
arrive quickly at our main results, we postpone the exact statements of 
these hypotheses for the next section. 
\begin{dfn} 
Let $f_1,\ldots,f_k\!\in\!\Z[x_1,\ldots,x_n]$, 
$F\!:=\!(f_1,\ldots,f_n)$,
and let $\sigma(F)$, the {\bf sparse size} (or {\bf bit size}) {\bf of $F$}, denote the 
number of bits\footnote{The sparse size thus obviously extends to 
integers: $a\!\in\!\Z \Longrightarrow \sigma(a)\!=\!1+\lceil 
\log_2 (1+|a|)\rceil$.} needed to write down the monomial term expansions of 
all the $f_i$.
(So only monomials which appear are counted, and we (internally) 
assume all coefficients {\bf and} exponents are written in binary.) \dia 
\end{dfn}
\begin{thm}
\label{thm:big} 
Following the notation above, suppose either GIPIT or DZH is true. 
Then $\hn\!\in\!\am$. Furthermore, 
\[ \text{GRH}\Longrightarrow\text{DZH}\Longrightarrow\text{GIPIT}
\Longrightarrow \hn\!\in\!\am \Longrightarrow 
[\hn\!\not\in\!\pp \Longrightarrow \pp\!\neq\!\np ],\] 
and both GIPIT and DZH can hold even if \grh is false. 
\end{thm} 

\noindent 
The very last implication follows easily from the basic 
facts that $\am\!\subseteq\!\pp^{\np^\np}$ and 
$\pp\!=\!\np \Longrightarrow \pp\!=\!\pp^\np=\!
\pp^{\np^\np}$ \cite{papa}. 
Note also that the implication $\hn\!\in\!\pp 
\Longrightarrow \pp\!=\!\np$ holds unconditionally, thanks to the aforementioned 
$\np$-hardness of $\hn$ \cite{plaisted}. 

The algorithm underlying Theorem \ref{thm:big} is based 
on the following quantitative result relating prime distributions and 
solvability over the complex numbers and finite fields. 
\begin{thm} 
\label{thm:stride} 
Following the notation above, there are positive integers 
$t_F$, $a_F$, and $C_F$ (depending on $F$)  
such that the following assertions hold: 
\begin{enumerate} 
\item{$F$ has no complex root $\Longrightarrow$ the 
mod $p$ reduction of $F$ has a root in $(\Z/p\Z)^n$ for at most 
$a_F$ primes $p$.  In particular, $\sigma(a_F)\!=\!O(\sigma(F)^2)$.}  
\item{$F$ has at least one complex root $\Longrightarrow$ 
for all integers 
$t\!\geq\!t_F$, there 
% C' instead of C (and remove some earlier %'s) if 
% the constants in Ono-Heath-Brown-Mitsui-Sokolovski\cup{i} 
% wind up being good enough! !!! 
is a prime $p\!\in\!\{t^{C_F}+1,\ldots,(t+1)^{C_F}\}$  
such that the mod $p$ reduction of $F$ has a root in $(\Z/p\Z)^n$.} 
\item{We have $C_F\!=\!2^{O(\sigma(F))}$ and $t_F\!=\!O(1)$ unconditionally.} 
\item{If either of DZH or GIPIT are true, then we can assert further that 
$C_F\!=\!O(\sigma(F)^{\kappa})$ and $t_F\!=\!O(\sigma(F)^{\kappa})$ 
for some effectively computable absolute constant $\kappa\!>\!0$. } 
\item{If GRH is true then we can assert even further that 
$C_F\!=\!O(\sigma(F)^{2.01})$ and $t_F\!=\!O(1)$.}   
\end{enumerate}
\end{thm} 

\noindent
{\bf Note:} {\em Throughout this paper, all $O$-constants 
and $\Omega$-constants are absolute and effectively computable. \dia }

\medskip 

Explicit formulae for $t_F$, $a_F$, and $C_F$, under the assumption of GRH,  
were derived earlier in \cite[Sec.\ 6.1.6, Rems.\ 12 and 13]{four}. Although  
Assertion (1) was observed earlier by Koiran \cite[Thm.\ 5]{hnam} (following 
easily from earlier work of J\'anos Koll\'ar \cite{kollar}), 
we point out that Assertion (2) appears be new. In particular, the existence 
of a sequence of polynomial growth, interleaving the set of 
primes $p$ for which the mod $p$ reduction of $F$ has a root in $(\Z/p\Z)^n$,  
was known earlier only under strong restrictions: One either 
needed to assume \grh \cite[Thm.\ 1]{hnam} or use unconditional results (e.g., 
those of \cite{heathbrown}) that would leave $t_F$ completely unknown and 
$C_F$ far too large to be algorithmically useful.  

Via a clever use of randomization, it is not hard to derive an $\am$ 
algorithm for $\hn$ directly from Theorem \ref{thm:stride}, assuming 
GIPIT or DZH. One may then wonder 
if one could dispense with randomness and instead simply check the solvability 
of $F$ over $\Z/p\Z$ for a number of primes polynomial in 
$\sigma(F)$, thus obtaining an algorithm in $\pp^\np$ for $\hn$. 
Unfortunately, the following family of examples shows that such a 
simplification is simply not possible.  
\begin{ex} 
\cite[Thm.\ 10]{hnam} 
For $n\!\geq\!1$, consider $F$ defined by the following two univariate 
equations: 
\[ x^{\pi_n}-1=0\]
\[ x-\pi_n=0 \] 
where $\pi_n$ is the product of the first $n$ prime numbers. 
Then $F$ clearly has no roots in $\C$ and $\sigma(F)\!=\!O(n\log n)$. 
(The latter fact follows easily from the Prime Number Theorem, 
stated in Section \ref{sub:ana} below.) However, the number of 
primes $p$ for which the mod $p$ reduction of $F$ has a root 
is at least $2^n$, which is clearly exponential in $\sigma(F)$. 
This abundance of primes follows from the fact (observed by 
Noam Elkies around 1996) that the number of 
prime factors of $\pi^{\pi_n}_n-1$ is at least $2^n$ \cite[Lem.\ 4]{hnam}. \dia 
\end{ex}  

\subsection{Riemann Hypotheses and Related Zeta Functions} 
\label{sub:ana}
Primordial versions of the connection between analysis 
and number theory are not hard to derive from scratch and have been known 
since the 19$\thth$ century. For example, letting 
$\zeta(s)\!:=\!\sum^\infty_{n=1} \frac{1}{n^s}$ 
denote the usual {\bf Riemann zeta function} 
(for any real number $s\!>\!1$), 
one can easily derive with a bit of calculus \cite[pp.\ 30--32]{des} that 
\[ \zeta(s)\!=\!\prod_{p \text{ prime}}\frac{1}{1-\frac{1}{p}} \text{ and 
thus } -\frac{\zeta'(s)}{\zeta(s)}\!=\!\sum^\infty_{n=1} 
\frac{\Lambda(n)}{n^s},\] 
where $\Lambda$ is the classical Mangoldt function 
which sends $n$ to $\log p$ or $0$, according as $n=p^k$ 
for some prime $p$ (and some positive integer $k$) or not. 
For a deeper connection, recall that $\pi(x)$ denotes the 
number of primes (in $\N$) $\leq\!x$ and that 
the {\bf Prime Number Theorem (PNT)} is 
the asymptotic formula $\pi(x)\sim\frac{x}{\log x}$ 
for $x\longrightarrow \infty$. 
Writing $\rho\!=\!\beta+i\gamma$ for real $\beta$ and $\gamma$, recall also 
that $\zeta$ admits an analytic continuation 
to the complex plane, sans the point $1$, via 
\[ \zeta(\rho)\!=2^\rho\pi^{\rho-1}\sin\left(\frac{\pi \rho}{2}\right)
\Gamma(1-\rho) \zeta(1-\rho),\]
where $\Gamma$ is Euler's famous continuous analogue of the factorial 
function \cite[Sec.\ 2]{des}. (We'll abuse notation henceforth by letting 
$\zeta$ denote the analytic continuation of $\zeta$ to 
$\C\setminus\{1\}$.) 
One should then note that the first proofs of PNT, 
by Hadamard and de la Vall\'ee-Poussin (independently, in 1896), were based 
essentially on the fact that $\zeta(\beta+i\gamma)$ has {\bf no} zeroes 
on the vertical line $\gamma\!=\!1$.\footnote{Shikau Ikehara 
later showed in 1931 that PNT is in fact {\bf equivalent} to the fact that 
$\zeta$ has no zeroes on the line $\gamma\!=\!1$ (the proof 
is reproduced in \cite{dym}).}   

Going to an even deeper connection, the {\bf Riemann Hypothesis (RH)} (from 
1859) is the following assertion: 

\begin{quote}
{\bf (RH)} The {\bf only} zeroes $\rho\!=\!\beta+i\gamma$ of $\zeta$ 
other than the ``trivial'' zeroes $\{-2,-4,-6,\ldots\}$ 
are those on the {\bf critical line} $\gamma\!=\!\frac{1}{2}$.
\end{quote}  
In particular, letting $\pi(x)$ denote the number of 
primes $p\!\leq\!x$, it is known that RH is true 
$\Longleftrightarrow \left|\pi(x)-\int^x_{2} \frac{dt}{\log t}\right|\!=\!
O(\sqrt{x}\log x)$ \cite{des}. 
Aside from a myriad of other hitherto unprovably sharp statements on the 
distribution of primes, the truth of RH is also known to imply a {\bf 
deterministic} polynomial-time algorithm for primality detection 
\cite{miller} even faster than the recent unconditional algorithm of 
Agrawal, Kayal, and Saxena \cite{aks}. 

Through the last estimate on 
$\pi(x)$, RH also easily implies that for any fixed $\eps\!>\!0$, there is a 
$C_\eps\!>0$ such that $\pi(x+x^{\frac{1}{2}+\eps})-\pi(x)\!\geq\!C_\eps\frac{x^{\frac{1}{2}
+\eps}} {\log x}$, for all $x$ sufficiently large. However, one can prove lower bounds nearly 
as good without RH. For example, it is known that 
$\pi(x+x^{0.535})-\pi(x)\!\geq\!\frac{x^{0.535}}{20\log x}$ for
all $x$ sufficiently large \cite{bakerharmann}. 
Preceding \cite{iwaniecjutila}, results of this flavor were derived via estimates on the sum 
\[S(x,T)\!:=\!\sum\limits_{\substack{\zeta(\rho)=0 \ , \ \rho=\beta+i\gamma\\ 0<\beta<1 \ 
 , \  |\gamma|<T}} \frac{x^{\rho}}{\rho}.\] 
In particular, since the work of Hoheisel \cite{hoheisel}, it has been known that 
sufficiently good bounds of the form $|S(x,T)|=o(x\log^2 T)$ yield results 
such as our last bound; {\bf independent of RH}. Note also that we at least know 
$|S(x,T)|\!=\!O(x\log^2 T)$ unconditionally \cite[Lem.\ 5.4 and Sec.\ 9]{lago} 
and that under RH we have $|S(x,T)|\!=\!O(\sqrt{x}\log^2 T)$. 

Let us now consider general number fields. 
\grh can be defined as follows: For any finite algebraic 
extension $K$ of $\Q$, define the {\bf Dedekind zeta 
function} via the formula $\zeta_K(s)\!:=\!\sum_{\ca} 
\frac{1}{(N\ca)^s}$, where $\ca$ ranges over all 
nonzero ideals of $\cO_K$, $N$ denotes the norm function, and $s\!>\!1$ \cite{bs}. (So 
$\zeta_\Q\!=\!\zeta$.) Then, like $\zeta$, the function $\zeta_K$ also 
admits an analytic continuation to $\C\!\setminus\!\{1\}$ \cite{lago}, and \grh is 
the following statement: 
\begin{quote}
{\bf (GRH)} The {\bf only} zeroes $\rho\!=\!\beta+i\gamma$ of $\zeta_K$ 
other than the ``trivial'' zeroes $\{-2,-4,-6,\ldots\}$ 
are those on the critical line $\gamma\!=\!\frac{1}{2}$.
\end{quote}  
In particular, letting $\pi_K(x)$ denote the number of {\bf prime} ideals 
$\cp$ of $\cO_K$ satisfying\\ 
\scalebox{.95}[1]{$N\cp\!\leq\!x$, 
it is known that \grh is true $\Longleftrightarrow 
\left|\pi_K(x)-\int^x_{2} \frac{dt}{\log t}\right|\!=\!O\left(\sqrt{x}
(n_K\log x+\log d_K)\right)$,}\\ 
where $n_K$ and $d_K$ respectively denote 
the field extension degree $[K:\Q]$ and absolute value of the discriminant of $K$ 
over $\Q$ (see, e.g., \cite[Thm.\ 1.1]{lago} and \cite{bs}). Also, 
in addition to numerous other important consequences in number theory and 
complexity theory, the truth of GRH implies sub-exponential speed-ups for detecting 
rational points on certain algebraic sets \cite{jcs}. 

Paralleling the special case $K\!=\!\Q$, GRH also easily implies that for any fixed 
\scalebox{.98}[1]{$\eps\!>\!0$, there is a $C_\eps\!>0$ such that 
$\pi_K\left(x+x^{\frac{1}{2}+\eps}\right)-\pi_K(x)\!\geq\!C_\eps\frac{x^{\frac{1}{2}+\eps}}
{\log x}$, for all $x\!>\!2$. (This}\\ follows routinely from, say, 
\cite[Thm.\ 1.1]{lago}.) However, trying to do without GRH 
introduces some difficulties: While it is known that one can always find 
positive constants $c_K$, $C_K$, and $t_K$ depending on $K$ such that 
$c_K\!<1/2$ and $\pi_K(x+x^{1-c_K})-\pi_K(x)\!\geq\!C_K\frac{x^{1-c_K}}{\log x}$ 
for all $x\!\geq\!t_K$ \cite{sokolovskii,mitsui,heathbrown,ono}, 
there are no known estimates for $C_K$ or $t_K$. In any event,  
in close analogy to the case of the ordinary prime counting function, 
bounds of this type can be derived by estimating the sum 
\[S_K(x,T)\!:=\!\sum\limits_{\substack{\zeta_K(\rho)=0 \ , \ \rho=\beta+i\gamma\\ 0<\beta<1 \
 , \  |\gamma|<T}} \frac{x^{\rho}}{\rho}.\]
Let us also note that $|S_K(x,T)|\!\leq\!O(x(\log d_K+n_K+\log T)\log T)$ 
unconditionally\\ 
\scalebox{.78}[1]{\cite[Lem.\ 5.4 and Sec.\ 9]{lago}, and that under GRH we have 
$|S_K(x,T)|\!\leq\!O(\sqrt{x}(\log d_K+n_K+\log T)\log T)$.} 

We are thus naturally led to DZH and GIPIT. 
These alternatives to GRH respectively deal with upper bounds on 
$|S_K(x,T)|$ and lower bounds on $\pi_K(x+y)-\pi_K(x)$.  

\bigskip 
\noindent
{\bf Dedekind Zero Hypothesis (DZH)}
% Wrote up about the same version on May 29, 2002 at the Courant 
% Institute, but I already had an uglier version earlier...
{\em Following the notation above, 
% let us define $N_K(b,T)$ to be the number of zeroes of $\zeta_K$ in the 
% strip $\{\beta+i\gamma \; | \; 
% b\!\leq\!\beta\!\leq\!1 \ , \ |\gamma |\!\leq\!T\}$. 
% Then there are absolute constants $A,c\!>\!0$ such that  
% the following quantitative statements on $\zeta_K$ hold: \\
% \medskip
% \noindent
% \fbox{$\stackrel{\text{{\bf Forbidden}}}{\text{{\scriptsize 
% {\bf Region}}}}$}  
% \scalebox{.86}[1]{$\left[ \gamma\!\geq\!(1+n_K+\log d_K)^c \text{ and } 
% \beta\!\geq\!1-\left\{ (2+\log(n_K\log d_K))^c \log^{1-A} \gamma \right\}^{-1} 
% \right]  \Longrightarrow \zeta_K(\beta+i\gamma)\!\neq\!0$.}\\
% \medskip
% \noindent
% \mbox{}\hspace{.2cm}\fbox{{\bf {\scriptsize Density}}}\hspace{.1cm}\mbox{} 
% \scalebox{.94}[1]{There is a nonnegative 
% function $\alpha : \R \longrightarrow \R$ with 
% $\lim\limits_{T\longrightarrow 
% +\infty} \frac{\alpha(T)}
% {\log^{A} T \left/ \left( 2+\log(n_K\log d_K) \right) ^c \right. }\!=\!0$}\\ 
% \mbox{}\hspace{1.9cm}\scalebox{.9}[1]{such that 
% $N_K(b,T)\!=\!O\left((1+n_K+\log d_K)^c
% T^{(1-b)(2+\log n_K)^c} \log^{\alpha(T)} T\right)$ for 
% $\frac{1}{2}\!\leq\!b\!\leq\!1$.}  }
% % Sarnak said that the lower bound on gamma in the forbidden region 
% % assertion may be the criticial issue... (around May 29, 2002) 
% %
% % note also that a portion of dzh can be derived 
% % as follows: \cite{heathbrown} has shown that 
% % (without unknown remaining constants) 
there exists an absolute constant $\bk\!>\!0$ such that 
$|S_K(x,T)|\!\leq\!x^{1-\frac{1}{1.01+(1+\log(n_K\log d_K))^\bk}}T^{\frac{1.99}
{1.01+(1+\log(n_K\log d_K))^\bk}}$ for all $x,T\!\geq\!(1+\log(n_K\log d_K))^\bk$. 
\dia } 

\bigskip 
\noindent
{\bf Good Short Interval Prime Ideal Theorem (GIPIT)}
{\em Following the notation above, there is an absolute constant $\kk\!>\!0$ such that\\ 
\scalebox{.82}[1]{$\displaystyle{\pi_K\left((x+1)^{1.01+(1+\log(n_K\log d_K))^\kk}\right)
-\pi_K\left(x^{1.01+(1+\log(n_K\log d_K))^\kk}\right)\!\geq\!1+\frac{n_K}{2}
(x+1)^{(1.01+(1+\log(n_K\log d_K))^\kk)/2}}$,}\\
for all $x\!\geq\!(1+\log(n_K\log d_K))^\kappa$. \dia } 

\noindent
Note that DZH and GIPIT become harder to prove as we decrease $\bk$ and 
$\kk$, and that DZH is the stronger assumption (cf.\ Theorem \ref{thm:big} of the last 
section). In particular, by our preceding discussion, the truth of GRH easily 
implies that both DZH and GIPIT are true for {\bf any} arbitrarily small positive $\bk$ and 
$\kk$. 
\begin{rem} 
Some experts have recently expressed the opinion that 
\grh maybe true modulo a small exception: the existence of so-called 
Landau-Siegel zeroes on the real line \cite{duke}. Note, however, that 
the existence of such zeroes does not affect the truth of DZH or GIPIT. \dia 
\end{rem}

\section{Background and Related Work } 
\label{sec:back} 
Let us first consider the following recent refined effective 
arithmetic version of Hilbert's Nullstellensatz.\\ 
{\bf Note:} {\em In what follows, we always normalize $n$-dimensional volume $\vol(\cdot)$ 
so that the $n$-dimensional unit cube has volume $n!$. \dia }  
\begin{thm} \cite[Thm.\ 4.1]{sombra} 
\label{thm:cool} 
Let $f_1,\ldots,f_k\!\in\!\Z[x_1,\ldots,x_n]$, $F\!:=\!(f_1,\ldots,f_k)$, 
$D\!:=\!\max_i\deg f_i$, and let $V_F$ denote the volume of the convex hull of\footnote{i.e., 
smallest convex set containing...} the union of 
$\{\bO,e_1,\ldots,e_n\}$ (the origin and the standard basis vectors of $\Rn$) and the set of 
all exponent vectors of the monomial term expansion of $F$. 
Then $F$ has {\bf no} roots in $\Cn \Longleftrightarrow$ 
there exist polynomials $g_1,\ldots,g_k\!\in\!\Z[x_1,\ldots,x_n]$ 
and a positive integer $a$ with $f_1g_1+\cdots+f_kg_k\!=\!a_F$ identically and
\begin{enumerate}
\item{$\deg g_i\!\leq\!2n^2 D V_F$} 
\item{$\sigma(a_F)\!\leq\!2(n+1)^3DV_F(\sigma(F)+\log k +14(n+1)D
\log (D+1))$ \qed}  
% \item{\scalebox{.85}[1]{$\sigma(g_1),\ldots,\sigma(g_k)\!\leq\!2m(n+1)^3DV_F(\log D)
% (\sigma(F)+\log k +14(n+1)D \log (D+1))$ \qed} } 
\end{enumerate} 
\end{thm} 

\noindent 
This result immediately implies that $\hn$ can be solved 
by solving a linear system with $O((2n^2DV_F)^n)$ 
variables and equations over the rationals (with total bit-size 
$O(\sigma(F)(2n^2DV_F)^n)$), thus 
easily yielding $\hn\!\in\!\expt$. That $\hn\!\in\!\pspa$ 
then follows immediately from the fact that linear algebra 
can be efficiently parallelized \cite{csanky}. 

Another result we'll need is an effective version of rational univariate reduction 
(a.k.a. effective primitive element theorem). 
\begin{thm}
\label{thm:unired} 
[(See \cite[Thms.\ 5 and 7]{four})]
Following the notation of Theorem \ref{thm:cool}, let\\ 
\scalebox{.9}[1]{$Z_F$ denote the 
zero set of $F$ in $\Cn$. Then there are univariate
polynomials $h_1,\ldots,h_n,\hat{h}_F\!\in\!\Z[t]$}\\ and positive integers 
$a_1,\ldots,a_n$ with the following properties:
\begin{enumerate}
\item{The number of irreducible components of $Z_F$ is bounded above
by $\deg \hat{h}_F$; and $\deg h_1,\ldots,
\deg h_n\!\leq\!  \deg \hat{h}_F\!\leq\!V_F$.}  
\item{For any root $\theta$ of $\hat{h}_F$ we have 
$F\left(\frac{h_1(\theta)}{a_1},\ldots,\frac{h_n(\theta)}{a_n}\right)\!=\!\bO$ and 
every irreducible component of $Z_F$ contains at least one point that can be expressed in 
this way.} 
\item{\scalebox{.85}[1]{$\sigma(a_1),\sigma(h_1)\ldots,\sigma(a_n),\sigma(h_n)\!=\!
O(V^5_F\sigma(\hat{h}_F))$ and $\sigma(\hat{h}_F)\!
=\!\cO(V_F[\sigma(F)+n\log D])$. \qed }}
% check this!!! 
\end{enumerate}
\end{thm}

\noindent
{\bf Sketch of Proof:} 
% Expand this later!!! 
The proof follows almost exactly the same development 
as \cite[Secs.\ 6.1.3 and 6.1.5]{jcs}, except that Theorem 23 there (an estimate on the 
maximum coefficient size of a variant of the toric resultant) is combined with  
\cite[Cor.\ 2.5]{sombra}. This allows us to omit a factor exponential 
in $n$ that would have otherwise appeared in Assertion (3) above. \qed

% some facts on complexity: 
% tarski (1930's/40's?): quantifier elim / R is decidable 
%                        (later shown dexp, and includes 
%                        computation of homology gps of 
%                        semi-alg sets...)
% 
% fischer & rabin (1974): presburger arithmetic (1929, 
%                         1st order theory of addition over N, 
%                         consistent and complete)  
%                         is doubly exponentially hard...  
%                         (note that divisibility and primes 
%                         are NOT expressible in presburger-ese...) 
%                      example thm: x<=y+1 -> y+2 > x
%
% mayr and meyer (1982): ideal membership in Q[x_1,...x_n] expspace-hard
% burgisser (1998): ...over any infinite field...

For technical reasons, it will be more convenient to work 
with a weighted variant of $\pi_K$ instead of $\pi_K$ itself. 
So let us define $\psi_K(x)\!:=\!\sum \log N\bp$ --- the natural generalization
of the Mangoldt $\Lambda$ function to prime ideals in algebraic
number rings --- where the sum ranges over all integral
ideals of the form $\bp\!=\!\bq^m$ with $\bq$ an unramified prime ideal
of $\cO_K$ and $N \bq^m\!\leq\!x$. The results from analytic number theory we'll 
need can then be summarized as follows.  
\begin{thm}
\label{thm:uwepit}
\scalebox{1}[1]{{\bf (Unconditional (Weighted) Effective Prime Ideal Theorem)}}
\mbox{}\\ 
For all $x,T\!\geq\!2$, we have\\ 
\scalebox{.78}[1]{$\displaystyle{\left|\psi_K(x)-\left(x-
\sum\limits_{\substack{\zeta_K(\rho)=0 \ , \ 
\rho=\beta+i\gamma\\ 0<\beta<1 \ , \  |\gamma|<T}} \frac{x^\rho}{\rho} +
\sum\limits_{|\rho|<\frac{1}{2}} \frac{1}{\rho}\right)\right|\!=\!
O\left(
\frac{x}{T}\left(n_K\log^2 x+\log d_K+n_K\log T\right)
+(n_K+\log d_K)\log x
\right)}$.}\\ 
In particular, for all $x\!\geq\!e^{10}\!\approx\!22027$ and $T\!\geq\!e^5\!\approx\!148.4$, 
we have that the left-hand side of the above asymptotic equality is no larger than\\ 
\scalebox{.68}[1]{$\displaystyle{\frac{x}{T}\left\{5n_K\log^2 x+36.55n_K\log x 
+375.2n_K\log T+196.6\log d_K + 351n_K\right\}
+(2\log d_K+3.51n_K)\log x+25\log d_K+283.1n_K}$. }  
\end{thm} 

\noindent
{\bf Sketch of Proof:} The first portion is nothing more than the special case of 
\cite[Thm.\ 7.1]{lago} where, in the notation there, $K\!=\!L$ and the 
conjugacy class $C$ of $\gal(L/K)$ is just the identity element $e$. (The 
sums over zeroes of $\zeta_K$ also simplify since there is just one underlying 
character of $H\!=\!\langle e \rangle$ and it is the trivial character.) The second 
portion is a variation of \cite[Thm.\ 9.1]{silberman}. (The latter theorem 
gives a bound on $\psi_K(x)-x$ instead, and assumes the truth of GRH.) In particular, 
by modifying the proof to {\bf avoid} the use of GRH (and employing 
\cite[Lem.\ 5.8]{silberman}), the result follows easily. \qed 

\noindent
We then immediately obtain the following corollary. 
\begin{cor}
\label{cor:uwesipit}
\scalebox{.82}[1]{{\bf (Unconditional (Weighted) Effective Short Interval 
Prime Ideal Theorem)}}\mbox{}\\ 
For all $x,y,T\!\geq\!2$, we have\\ 
\scalebox{.63}[1]{$\displaystyle{\psi_K(x+y)-\psi_K(x)\!=\!y-  
\sum\limits_{\substack{\zeta_K(\rho)=0 \ , \ 
\rho=\beta+i\gamma\\ 0<\beta<1 \ , \  |\gamma|<T}} \frac{(x+y)^\rho-x^\rho}{\rho} - 
O\left(\frac{x+y}{T}\left(n_K\log^2(x+y)+\log d_K+n_K\log T\right)
+(n_K+\log d_K)\log(x+y)
\right)}$.}\\ 
In particular, for all $x,y\!\geq\!e^{10}\!\approx\!22027$ and 
$T\!\geq\!e^5\!\approx\!148.4$, we have that the equality sign above can be replaced by 
a $\leq$ sign and the $O(\cdot)$ estimate replaced by the following:\\ 
\scalebox{.55}[1]{$\displaystyle{\frac{x+y}{T}\left\{10n_K\log^2(x+y)+73.1n_K\log(x+y) 
+750.4n_K\log T+393.2\log d_K + 702n_K\right\}
+(4\log d_K+7.02n_K)\log(x+y)+50\log d_K+566.2n_K}$. \qed }  
\end{cor} 

While we won't need the following estimate to prove our main results, we 
include it so the reader can see how the conditional estimates on 
$\pi(x+y)-\pi(x)$ and $\pi_K(x+y)-\pi_K(x)$ from Section 
\ref{sub:ana} follow immediately from our unconditional corollary above, should 
GRH turn out to be true. 
\begin{lemma} 
(See \cite[Lem.\ 5.4]{lago} and \cite[Lem.\ 5.7 and the Proof of 
Thm.\ 9.1]{silberman}) 
We have $\displaystyle{\sum\limits_{\substack{\zeta_K(\rho)=0 \ , \
\rho=\beta+i\gamma\\ 0<\beta<1 \ , \  |\gamma|<T}} \frac{1}{|\rho|}\!=\!
O(\log^2 T + (n_K+\log d_K)\log T)}$. In particular, for all $T\!\geq\!e^5\!\approx\!148.4$, 
we have 
\scalebox{1}[1]{$\displaystyle{\sum\limits_{\substack{\zeta_K(\rho)=0 \ , \
\rho=\beta+i\gamma\\ 0<\beta<1 \ , \  |\gamma|<T}} \frac{1}{|\rho|}\!\leq\!
3.1\log^2 T + (77.1n_K+8\log d_K)\log T }$. \qed}
\end{lemma} 

The following basic result will help us relate $\log d_K$ and $n_K$ to the 
bit-size of the defining polynomial of $K$. Recall that 
the square-free part of a univariate polynomial $f$ is 
$f/\gcd(f,f')$. 
\begin{prop} 
\label{prop:disc}  
(See \cite[Lemmata 2.1 and 4.1]{jcs}) 
Suppose $K\!=\!\Q[x_1]/\langle f \rangle$ where $f$ is the square-free part of some polynomial 
$g\!\in\!\Z[x_1]\!\setminus\!\Z$. Then 
\[n_K\!\leq\!\deg g \text{ \ and \ } \log d_K\!=\!O((\deg g)
\sigma(g)+(\deg f)^2).\] More precisely,\\ 
\scalebox{.72}[1]{$\displaystyle{\log d_K\leq (2\deg g -1)(\sigma(g)+(\deg g 
+\alpha)\log 2)+\frac{2\deg g-1}{2}\log (\deg g+1)+\frac{\deg g}{2}\log((\deg g)
(2\deg g+1)/6)}$. \qed}
\end{prop} 

Finally, let us briefly and informally review the complexity classes mentioned earlier 
in the introduction. First, our underlying computational model will be the 
classical Turing machine, which for our purposes can be assumed to be a lap-top 
with infinite memory and a flawless operating system. In particular, bit operations 
(and time) can be identified with the number of clock-cycles used by the underlying 
microprocessor. Input size, for an input such as $F$, is then simply $\sigma(F)$, 
and can also be interpreted as the amount of memory used by our idealized lap-top to 
store $F$.

Our complexity classes can then be summarized as follows. Note that 
an {\bf oracle in $\pmb{A}$} is a special machine which is allowed to run 
an algorithm with complexity $A$ in one unit of time, and all problems below 
are {\bf decision} problems. 
\begin{itemize} 
\item[$\pp$]{ The family of problems which can be done within time 
polynomial in the input size.\footnote{Note that the underlying polynomial 
depends only on the {\bf problem}, not the particular instance of the 
problem.}} 
\item[$\rp$]{ The family of problems admitting a polynomial-time algorithm 
for which a {\tt ``Yes''} answer is always correct but a {\tt ``No''} answer 
is wrong with probability $\frac{1}{2}$. } 
\item[$\corp$]{ The family of problems admitting a polynomial-time algorithm 
for which a {\tt ``No''} answer is always correct but a {\tt ``Yes''} answer 
is wrong with probability $\frac{1}{2}$. } 
\item[$\np$]{ The family of problems where a {\tt ``Yes''} answer can be 
{\bf certified} within time polynomial in the input size.} 
\item[$\conp$]{ The family of problems where a {\tt ``No''} answer can be 
{\bf certified} within time polynomial in the input size.} 
\item[$\am$]{ The family of problems solvable by a $\bpp$ algorithm which has been 
augmented with exactly {\bf one} use of an oracle in $\np$. } 
\item[$\np^\np$]{ The family of problems where a {\tt ``Yes''} answer can be 
certified by using an $\np$-oracle a number of times polynomial in the 
input size. } 
\item[$\conp^\np$]{ The family of problems where a {\tt ``No''} answer can be 
certified by using an $\np$-oracle a number of times polynomial in the 
input size. } 
\item[$\pp^{\np^\np}$]{ The family of problems solvable within time polynomial in the 
input size, with as many calls to an $\np^\np$ oracle as allowed by the time bound. } 
\item[$\pspa$]{ The family of problems solvable within time polynomial in the 
input size, provided a number of processors exponential in the input size is allowed. } 
\item[$\expt$]{ The family of problems solvable within time exponential in the 
input size.} 
\end{itemize} 

We note that $\rp$ algorithms and $\bpp$ algorithms are sometimes respectively referred to 
as {\bf Monte-Carlo} and {\bf Atlantic City} algorithms. 

\section{Reducing Theorem \ref{thm:big} to Theorem \ref{thm:stride}} 
\label{sec:big} 

Let us begin with a description of the algorithm 
which proves Theorem \ref{thm:big}:\\
 
\noindent
{\sc Algorithm} {\tt KAMHN}\footnote{Short for Koiran-Arthur-Merlin-Hilbert-Nullstellensatz. 
A similar algorithm was outlined on Page 6 of \cite{hnam} (but {\bf not} in the 
journal version of the paper), hence our naming.}  
{\em 
\begin{itemize}
\item[{\bf Input}]{A system $F$ of $k$ polynomials 
in $n$ variables with integer coefficients} 
\item[{\bf Output}]{A declaration, correct with probability 
$\geq\!\frac{2}{3}$, of whether $F$ has a complex root 
or not.} 
\item[{\bf Step 0}]{Let $t_F$, $a_F$, and $C_F$ be 
the integer constants from Theorem \ref{thm:stride}. } 
\item[{\bf Step 1}]{ Pick a (uniformly distributed) random 
integer $t\!\in\!\{t_F,\ldots,t_F+3a_F\}$.} 
\item[{\bf Step 2}]{ \scalebox{.91}[1]{Using an $\np$ oracle just once, 
decide if there is a prime $p\!\in\!\{t^{C_F}+1,\ldots,  
(t+1)^{C_F}\}$}\\ 
such that the mod $p$ reduction of $F$ has a root in 
$(\Z/p\Z)^n$. If so, declare that $F$ has a complex root. Otherwise, 
declare that $F$ has no complex root. \dia  } 
\end{itemize} } 
% a ``NO'' answer is always correct, i.e., 
% algor always works if F actually has a root...
% otherwise, fails with prob 1/3...

The last technical lemma we will need relates short interval bounds 
between $\pi_K$ and $\psi_K$. 
\begin{lemma} 
\label{lemma:red} 
Suppose $f\!\in\!\Z[x_1]\!\setminus\!\Z$ is square-free and 
$K\!:=\!\Q[x_1]/\langle f\rangle$. Then, assuming $x\!\geq\!1+\log d_K$ and $y\!\geq\!0$,\\ 
\scalebox{.86}[1]{$\displaystyle{\psi_K(x+y)-\psi_K(x)\!\geq\!\left(1+\frac{7n_K}{2}\sqrt{x+y}
\right)\log(x+y) \Longrightarrow \pi_K(x+y)-\pi_K(x+y)\!\geq\!1+\frac{n_K}{2}\sqrt{x+y}}$.} 
\end{lemma} 

\noindent
{\bf Proof:} Note that the right-hand bound can be enforced via a 
similar bound on yet another related weighted prime power counting function: 
First, define 
$\Theta_K(x)\!:=\!\sum \log N\bp$ where the sum ranges over all 
unramified prime ideals of $\cO_K$ of norm $\leq\!x$. Since 
\[\Theta_K(x+y)-\Theta_K(x)\!\leq\!(\pi_K(x+y)-\pi_K(x))\log(x+y),\] 
the right-hand bound of our lemma is clearly implied by 
\[ \mbox{}\hspace{-.9in}(\heartsuit)\hspace{.9in} 
\Theta_K(x+y)-\Theta_K(x)\!\geq\!\left(1+\frac{n_K\sqrt{x+y}}{2}\right)
\log(x+y).\] 
Now, 
\[\sum\limits_{\substack{(\bp,m)\\ m\geq 2 \\ N\bp^m \leq x} } 
\log N\bp = O(\sqrt{x})\leq \sum\limits^{\lfloor \log_2 x\rfloor}_{k=2} n_K x^{1/k}\log x 
\leq n_k(\log x)\sum\limits^{\lfloor \log_2 x\rfloor}_{k=2} x^{1/k}\]
\[\leq n_k(\log x)(\sqrt{x}+x^{1/3}\log_2 x)\leq n_K(\log x)(3\sqrt{x})=3n_K\sqrt{x}\log x.\] 
So we then clearly have $0\!\leq\!\psi_K(x)-\Theta_K(x)\!=\!3n_K\sqrt{x}\log x$, 
and it thus suffices to enforce 
\[ 
\psi_K(x+y)-\psi_K(x)\geq \left(1+\frac{7n_K\sqrt{x+y}}{2}\right)
\log(x+y).\]
So we are done. \qed 

\noindent
{\bf Proof of Theorem \ref{thm:big}}: 
Clearly, to prove the first assertion of Theorem \ref{thm:big}, it suffices to prove 
that the truth of either DZH or GIPIT implies that the algorithm above is correct and 
runs within $\am$. We will handle the remaining implications later. 

To see correctness, simply note that Theorem \ref{thm:stride} implies that 
$F$ has no complex roots (resp.\ has at least one complex root) $\Longrightarrow$ 
at most one third (resp.\ {\bf all}) of the half-open 
intervals $\{[t^{C_F}+1,(t+1)^{C_F})\}_{t \in 
\{t_F,\ldots,t_F+3a_F\} }$ has (resp.\ have) a prime $p$ such that 
the mod $p$ reduction of $F$ has a root in $(\Z/p\Z)^n$. 
So a declaration of {\tt ``No Complex Roots''} is always 
correct, while a declaration of {\tt ``There is at Least 
One Complex Root''} is wrong with probability $\frac{1}{3}$. 
So correctness is actually independent of DZH or GIPIT or any unproven assumption. 

To see that {\tt KAMHN} runs in $\am$, let us now assume either DZH or GIPIT.  
Note then that by Assertion (2) of Theorem \ref{thm:stride}, no integer computed 
has bit size larger than $O(\sigma(F)^{\kappa+2})$, so 
a single use of an $\np$ oracle really does suffice to run our algorithm. 
Note also that since $\kappa$ is effectively computable, the constants 
$t_F$, $a_F$, and $C_F$ are all computable in polynomial time by recursive 
squaring. So {\tt KAMHN} indeed runs in $\am$. 

We must now prove the remaining assertions. The implications 
$[$GRH$\Longrightarrow$DZH$]$ and 
$[\hn\!\in\!\am \Longrightarrow [\hn\!\not\in\!\pp \Longrightarrow \pp\!\neq\!\np ]]$ 
were already proved during our discussion in Sections 1 and \ref{sub:ana}. 
The same goes for the fact that GIPIT and DZH can still hold even if GRH is false. 
So, thanks to our earlier application of Theorem \ref{thm:stride}, we need only prove 
DZH$\Longrightarrow$GIPIT. 

To prove the last implication, note that DZH implies that 
\[|S(x^C,x^{C/2})|\!\leq\!x^{C\left(1-\frac{0.005}{1.01+(1+\log(n_K\log d_K))^\bk}\right)}, 
\text{ for all } x\!\geq\!(1+\log(n_K\log d_K)^{\bk}.\] 
So, provided $C\!>\!200(1.01+(1+\log(n_K\log d_K))^\bk)$, we obtain that 
$|S(x^C,x^{C/2})|\!\leq\!x^{C-1}$ for all $x\!\geq\!(1+\log(n_K\log d_K)^{\bk}$. 
Let us define $C$ a bit larger asymptotically by 
setting $C\!:=\!(1.01+(1+\log(n_K\log d_K))^\kk$ for any $\kk\!>\!\bk$.

Now recall that $Cx^{C-1}\!\leq\!(x+1)^C-x^C\!\leq\!C(x+1)^{C-1}$. 
Corollary \ref{cor:uwesipit}, with $T\!:=\!x^{C/2}$, then yields\\ 
\scalebox{.75}[1]{$\displaystyle{\psi_K\left((x+1)^C\right)-\psi_K(x^C)\!=\!
C(x+1)^{C-1}-o(x^{C-1})-
O(x^{1-\frac{C}{2}}(Cn_K\log^2(x+1)+n_K\log x+C(\log d_K)\log x))}$,}\\ 
\[ = \Omega(x^{C-1}) \]
where the underlying constants are effectively computable and the implied 
inequality takes effect for all $x\!\geq\!(1+\log(n_K\log d_K))^{\kappa'}$ for 
some absolute effectively computable $\kappa'\!>\!0$. In particular, 
by construction $C-1\!>\!C/2$ and we thus obtain 
\[\psi_K\left((x+1)^C\right)-\psi_K(x^C)\!\geq\!\left(1+\frac{7n_K}{2}(x+1)^{C/2}\right)\log(x+1)\]
for all $x\!\geq\!(1+\log(n_K\log d_K))^{\kappa''}$, for
some absolute effectively computable $\kappa''\!>\!0$. Replacing $\kappa$ by 
$\max\{\kappa,\kappa''\}$, and substituting $x+y\leftarrow (x+1)^C$ and 
$x\leftarrow x^C$ into Lemma \ref{lemma:red}, we are done. \qed 

\section{Prime Distribution and Proving Theorem \ref{thm:stride}}  
\label{sec:stride} 

\noindent
{\bf Proof of Theorem \ref{thm:stride}:} Assertion (1) 
follows easily from the following argument which has 
appeared earlier in various incarnations, e.g., \cite{hnam,hmps,four}: 
By assumption, Theorem \ref{thm:cool} tells us that the mod $p$ reduction 
of $F$ has a root in $(\Z/p\Z)^n \Longrightarrow p$ divides $a_F$. Since
the number of prime divisors of an integer $a$ is no more than
$1+\log a$ (since any prime power other than $2$ is bounded below by
$e\!<\!2.718281829$), an elementary calculation yields Assertion (1).  

The proof of Assertions (2) relies on effective univariate reduction 
and some density estimates for prime ideals. In particular, 
we will derive Assertion (2) by examining the following cleaner abstraction: 
Let $\pi_F(x)$ denote the number of primes $p\!\leq\!x$ such that 
the mod $p$ reduction of $F$ has a root in $(\Z/p\Z)^n$. 
Clearly then, it is enough to find sufficiently good upper bounds on $y$ (as 
a function of $x$) and $t_F$ such that for all $x\!\geq\!t_F$,  
\[ \mbox{}\hspace{-1.75in}(\star) \hspace{1.75in}\pi_F(x+y)-\pi_F(x)\!\geq\!1. 
\] 
To do this, we will pass through various 
reductions, eventually culminating in enforcing a 
particular lower bound on $\pi_K(x+y)-\pi_K(x)$, for some number field $K$ depending 
naturally on $F$. In particular, GIPIT will immediately imply a lower bound 
good enough to imply Assertion (3), and without GIPIT we can still derive the 
weaker Assertion (4) by employing older theorems on the distribution of prime ideals 
in $O_K$. (Since Theorem \ref{thm:big} contains the implication 
DZH$\Longrightarrow$GIPIT, we need not worry about assuming DZH any more.)  
Our framework will then imply Assertion (5) almost trivially. 

To start, let $h_1,\ldots,h_n,\hat{h}_F\!\in\!\Z[t]$ and 
$a_1,\ldots,a_n$ be the polynomials and integers arising from 
applying Theorem \ref{thm:unired} to our input $F$. 
Letting $f$ be the square-free part of $\hat{h}_F$, note that
$p\!\not\!| \; \mathrm{lcm}(a_1,\ldots,a_n)$ {\bf and}
the mod $p$ reduction of $f$ has a root in $(\Z/p\Z)^n \Longrightarrow$ 
the mod $p$ reduction of $F$ has a
root in $(\Z/p\Z)^n$. So to prove ($\star$) it clearly 
suffices to instead enforce 
\[ \mbox{}\hspace{-.7in}(\star\star) \hspace{.7in}
N_f(x+y)-N_f(x)\!\geq\!1, \text{ for all } x\!\leq\!1+\log d_K\] 
where $N_f(x)$ denotes the number of linear factors
of the mod $p$ reduction of $f$ {\bf summed over all primes
$p\!\leq\!x$}. In particular, since  
$\sigma(f)\!=\!\sigma(\hat{h}_F)+(V_F+c_0)\log 2$, where 
$c_0\!=\!2-\frac{3}{4\log 2} \!\leq\!0.91798$ \cite[Lemma 5]{four}, 
Theorem \ref{thm:unired} and Proposition \ref{prop:disc} tell us that 
Assertion (2) is true if we can show that we can take $y\!=\!x^{1-1/C_F}$ 
in $(\star)$. However, the latter condition was already noted as true in our 
discussion of unconditional lower bounds on $\pi(x+y)-\pi(x)$ from Section 
\ref{sub:ana}. Similarly, our discussion of GRH-dependent lower bounds 
immediately implies Assertion (5). Assertion (3) then routinely from \cite[Thm.\ 1.3]{lago}, 
which contains an unconditional version of the Prime Ideal Theorem. 

To prove Assertion (4) we must now show GIPIT implies that we can find $t_F$ and $C_F$ with 
$\sigma(t_F),\sigma(C_F)\!=\!O(\sigma(f)^\kappa)$ for some  
effectively computable absolute constant $\kappa\!>\!0$. 
Letting $K\!:=\!\Q[x_1]/\langle f\rangle$, note now that $p\!\not\!| \; d_K$ implies that 
the mod $p$ reduction of $f$ has a linear factor 
iff $p$ lies below a prime in $\cO_K$ of residue class degree $1$ 
\cite[Ch.\ I, Sec.\ 8, Prop.\ 25, pp.\ 27--28]{lang}. Note also that $\bp$ (of residue class 
degree $1$) lies over $p \Longrightarrow N\bp\!=\!p$. Now let 
$\pi^1_K(x)$ denote the number of 
unramified primes $\bp$ of $\cO_K$ with $N\bp\!\leq\!x$ and residue class degree $1$. 
It is then clear that $\left|\pi^1_K(x)-N_f(x)\right|\!\leq\!(n_K-1)
(1+\log d_K)$ (since the correspondence between factors and primes breaks down 
only for $p|\Delta_K$) and the left-hand side is {\bf constant} 
for all $x$ larger than any prime dividing $d_K$.  
So we can clearly replace ($\star\star$) by the combined assertions 
\[ \mbox{}\hspace{-1.75in}(\clubsuit)\hspace{1.75in} 
\pi^1_K(x+y)-\pi^1_K(x)\!\geq\!1\] 
and 
\[ \mbox{}\hspace{-1in}(\clubsuit\clubsuit)\hspace{1in} \text{ no prime } \geq\!t_F 
\text{ divides any } a_1,\ldots,a_n \text{ {\bf or }} \Delta_f.\] 

We are now ready for our final reduction: Let 
$\pi_K(x)$ denote the number of unramified primes $\bp$ of 
$\cO_K$ with $N\bp\!\leq\!x$. Note then that 
for any prime $\bp$ of $\cO_K$ of residue class field 
degree $\geq\!2$ lying over $p$, we must have $N\bp\!\geq\!
p^2$, which in turn implies $p\!\leq\!\sqrt{x}$. Since the 
number of such primes $\bp$ over $p$ is $\leq\!n_K/2$, we must 
then have $0\!\leq\!\pi_K(x)-\pi^1_K(x)\leq\!\frac{n_K\sqrt{x}}{2}$. 
So it suffices to replace ($\clubsuit$) by 
\[ \mbox{}\hspace{-1.4in}(\spadesuit)\hspace{1.4in} \pi_K(x+y)-
\pi_K(x)\!\geq\!1+\frac{n_K\sqrt{x+y}}{2}.\] 

Via GIPIT, we can then substitute $x+y\leftarrow (x+1)^{C_F}$ and
$x\leftarrow x^{C_F}$ to immediately obtain $(\spadesuit)$. 
Taking $t_F$ to be the smallest power of $2$ larger than $d_K$ 
(which, by Proposition \ref{prop:disc} and Theorem \ref{thm:unired}, 
must have size $O(1+\log(n_K\log d_K))^2)$), we can then immediately 
obtain $(\clubsuit\clubsuit)$. Note in particular that 
$t_F$ can be computed in polynomial time by recursive squaring. 
So we are done. \qed 

\section*{Acknowledgements} 
Around March of 2002, Enrico Bombieri kindly took the time to look at one  
of my earlier papers (\cite{four}) and spotted an error in 
one of the proofs. (The result in question was Part (b) of 
Theorem 4 and the error changes the numerical estimates 
of Remark 13 in Section 6.1.6 but does not invalidate the 
result.) The present paper arose when I checked the details and 
noted that the earlier assumptions of GRH could be considerably 
weakened. I am therefore deeply indebted to Professor Bombieri 
for helping inspire the present paper. 
Later, around May 2002, I had the wonderful good fortune of meeting Lior 
Silberman at the Workshop on Zeta-Functions and Associated Riemann 
Hypotheses at the Courant Institute of Mathematical Sciences, New York University. 
Lior provided a copy of \cite{silberman} just when I needed some extremely 
explicit $\check{\text{C}}$ebotarev-like estimates. I am thus also indebted to 
Professor Silberman for writing and sharing his notes with such 
generosity and good timing. 

\bibliographystyle{acm}

\end{document}